\documentclass[11pt]{article}
\usepackage{geometry}
\usepackage{indentfirst}
\usepackage{algorithm}
\usepackage{algorithmic}
\usepackage{amsthm}
\usepackage{amsmath}
\usepackage{amstext}
\usepackage{amssymb}
\usepackage{dsfont}
\usepackage{epstopdf}
\usepackage{subcaption}
\usepackage{booktabs}
\usepackage{cite, bm}
\usepackage[colorlinks,linkcolor=blue,citecolor=blue,anchorcolor=blue]{hyperref}
\usepackage[braket]{qcircuit}
\usepackage{braket}
\newtheorem{theorem}{Theorem}
\newtheorem{lemma}{Lemma}
\newtheorem{corollary}{Corollary}
\newtheorem{definition}{Definition}

\title{Generalized quantum singular value transformation with application in quantum conjugate gradient least squares algorithm
}

\author{Yu-Qiu Liu\footnotemark[2] \footnotemark[1]
	\and Hefeng Wang\footnotemark[3] \footnotemark[1]
    \and Hua Xiang\footnotemark[2] \footnotemark[4] \footnotemark[1]
    }

\date{}

\begin{document}
	\maketitle
	\renewcommand{\thefootnote}{\fnsymbol{footnote}}
	\footnotetext[2]{School of Mathematics and Statistics, Wuhan University, Wuhan 430072, China.}
    \footnotetext[3]{Department of Applied Physics, School of Science, Xi’an Jiaotong University and Shaanxi Province Key Laboratory of Quantum Information and Quantum Optoelectronic Devices, Xi’an, 710049, China}
    \footnotetext[4]{Hubei Key Laboratory of Computational Science, Wuhan University, Wuhan 430072, China}
	\footnotetext[1]{E-mail address: {\tt 15516557320@163.com (Y. Q. Liu ), } {\tt wanghf@mail.xjtu.edu.cn (H. Wang), }{\tt hxiang@whu.edu.cn (H. Xiang)}.}
\begin{abstract}
Quantum signal processing (QSP) and generalized quantum signal processing (GQSP) are essential tools for implementing the block encoding of matrix functions. The achievable polynomials of QSP have restrictions on parity, while GQSP eliminates these restrictions. But GQSP only constructs functions of unitary matrices. In this paper, we further investigate GQSP and extend it to general matrices. Compared with the quantum singular value transformation (QSVT), our proposed method relaxes the requirements on the parity of polynomials. We refer to this extension as generalized quantum singular value transformation (GQSVT). Subsequently, by utilizing the relationship between generalized matrix functions and standard matrix functions, we propose a classical-quantum hybrid quantum conjugate gradient least squares (CGLS) algorithm using GQSVT.   

\end{abstract}

\section{Introduction}
QSP is one of the most efficient and versatile methods for approximate matrix functions in recent years, requiring only a few ancillary qubits. QSP and its extensions have become important tools in quantum computing, with applications including ground-state preparation \cite{YLY2022Ground, DL2024MultilevelQSP, MGB2025Coherent}, solving linear systems \cite{LT20filtering, TWY2024QCG}, Hamiltonian simulation \cite{LC2017Uniform, LC2017Optimal, LC19Hamiltonian}, and block encoding \cite {CLV2024Explicit, LNY2023Efficient, SCC2024Block-encoding}.  It was originally proposed by Low et al. in 2016 to design composite quantum gates that can implement quantum response functions \cite{LYC16Methodology}. These functions are polynomials with definite parity and satisfy the condition that magnitudes do not exceed $1$. In 2019, Low and Chuang presented a technique known as qubitization \cite{LC19Hamiltonian}. This method embeds a Hamiltonian $H$ in an invariant subspace and can compute operator functions of $H$ with optimal query complexity. Utilizing qubitization and singular value decomposition (SVD), QSVT \cite{GSL19QSVT} further extends QSP to embed the polynomial of general operators into a unitary. Subsequently, GQSP was proposed in 2024 by Motlagh and Wiebe \cite{MW24GQSP}. It removes the restrictions on the parity of implementable polynomials in QSP, retaining only the condition that the polynomial $p$ satisfies $|p|\leq 1$. This relaxation further broadens the applicability of QSP.

However, GQSP only constructs the function of unitary $U=e^{iH}$, where $H$ is a Hermitian matrix. In this paper, we generalize GQSP from unitary matrices to general matrices by combining singular value decomposition and qubitization techniques. We adopt the framework of GQSP, using arbitrary $SU(2)$ rotations as the signal processing operator and controlled unitary matrices as the signal operator. We use two classes of $0$-controlled and $1$-controlled unitary operators to implement transformations between different spaces. Using the Chebyshev expansion of arbitrary polynomials, our method further extends the range of realizable block encoding matrix functions, that is, the block encoding of polynomials without definite parity of general operators.

QSP and its series of extensions have provided new and powerful approaches for solving the quantum linear system problem (QLSP). Solving linear systems is among the most frequent problems in scientific computing. The classical methods for solving linear systems are classified as either direct or iterative \cite{Vorst2003, Golub2013Matrix, Meurant2020Nonsymmetric, Saad2003, Barrett1994Templates}. However, as the scale of the problem gradually increases, the computational resources and time required by classical computers to solve such linear systems increase greatly. Quantum computers hold the promise of enabling new algorithms that can solve problems requiring excessive resources on classical computers \cite{Xiang2022, Nielsen2010, Shor97Polynomial, HHL09, Grover96fast}. 

The QLSP refers to the preparation of a quantum state $\ket{\bm{x}}=A^{-1}\ket{\bm{b}}$, which is proportional to the solution of a given linear system $A\bm{x}=\bm{b}$. Harrow, Hassidim, and Lloyd presented the HHL algorithm \cite{HHL09} in 2009, which was the first quantum algorithm to solve linear systems. The computational complexity of the HHL algorithm is $O(\kappa^2\log(n)/\varepsilon)$ and scales logarithmically with the system size $n$, where $\kappa=\|A\|\|A^{-1}\|$ is the condition number of $A$, $\|\cdot\|$ is the usual $\ell^2$  matrix norm, and $\varepsilon$ is the desired precision. In subsequent studies, the HHL algorithm has been further improved \cite{AA10VTAA, CKS17LCU, WZP18Dense}. The development of QSP has introduced new methods for solving QLSP. Gily\'{e}n et al. found the polynomial that approximates $1/x$ and employed QSVT to construct an operator that approximates $A^{-1}$ with maximum circuit depth $O(\kappa \log(\kappa/\varepsilon))$\cite{GSL19QSVT, MRT21Grand}. Lin and Tong presented a quantum eigenstate filtering algorithm and used it to solve quantum linear system problems, achieving near-optimal $O(\kappa\log(1/\varepsilon))$ complexity for both $\kappa$ and $\varepsilon$ \cite{LT20filtering}. In 2024, Toyoizumi et al. introduced a quantum conjugate gradient method that achieves a square-root improvement for $\kappa$ in the maximum circuit depth\cite{TWY2024QCG}.

As an application of GQSVT, we use it to solve the least squares (LS) problem, a fundamental and widely used approach for data fitting\cite{Trefethen97Numerical}. For a matrix $A \in \mathbb{R}^{m \times n}$ and a vector $\bm{b}\in \mathbb{R}^{n}$ $(m\geq n)$, the LS procedure find the vector $\bm{x}$ that minimizes $\|A\bm{x}-\bm{b}\|$. The vector $\bm{x}$ is the LS solution if and only if $A^{\top}A\bm{x} = A^{\top}\bm{b}$ when $A$ has full column rank. This equation can be solved by the conjugate gradient (CG) method, which leads to the CGLS algorithm. The CG method works for symmetric positive (SPD) definite problems and it is an iterative Krylov subspace method, where the Krylov subspace with dimension $j$ is defined as $\mathcal{K}_j(A,\bm{v})=\operatorname{span}\{\bm{v}, A\bm{v}, A^2\bm{v}, \ldots, A^{j-1}\bm{v}\}$. 

This paper is organized as follows. In Sect.2, we present several fundamental definitions and a brief review of the GQSP method. In Sect.3, we provide a detailed process for extending GQSP to GQSVT. We also discuss its circuit implementation. In Sect.4, we first introduce the LS problem and the CGLS algorithm, then we present the quantum CGLS algorithm by employing GQSVT to compute vectors and using the swap test to calculate inner products. By using the convergence properties of CG, we derive the maximum circuit depth for our proposed algorithm. Conclusions are given in Sect.5.

\section{Preliminaries}
\subsection{Functions of matrices}
We first review the definitions of standard matrix functions and generalized matrix functions.

\textbf{Standard matrix functions.} For simplicity, let $A$ be an $n \times n$ Hermitian matrix with the eigendecomposition $A = U\Lambda U^{\dagger}$, where $U$ is unitary and $\Lambda = \operatorname{diag} (\lambda_1, \lambda_2, ..., \lambda_n)$ is a diagonal matrix containing the eigenvalues of $A$. For a function $f$, the standard matrix function $f(A)$ is defined as 
$$f(A) = U f(\Lambda) U^{\dagger}, $$ 
where $f(\Lambda) =\mathrm{diag} (f(\lambda_1), f(\lambda_2), ..., f(\lambda_n))$.

\textbf{Generalized matrix functions.} Let $A$ be an $m \times n$ matrix with rank $r = \operatorname{rank}(A)$. The singular value decomposition of $A$ is given by $A=W\Sigma V^{\dag}$, where $W$ and $V$ are unitary. Here, $\Sigma= \mathrm{diag} (\sigma_1, ..., \sigma_r, 0, ..., 0)$ is an $m \times n$ diagonal matrix whose diagonal entries $\sigma_k$ are the non-negative singular values of $A$. For a function $f$, the generalized matrix functions $f^{\diamond}(A)$ is defined as
$$f^{\diamond}(A)=Wf(\Sigma)V^{\dagger},$$
where $f(\Sigma)= \operatorname{diag} (f(\sigma_1), ..., f(\sigma_r), 0, ..., 0)$.
The left and right generalized matrix polynomials are defined as
$$f^{L}(A)=Wf(\Sigma)W^{\dagger}~~\text{and}~~f^{R}(A)=Vf(\Sigma)V^{\dagger},$$
respectively\cite{HB73generalizedMF, BG2003Generalized}.

In general, $f^{\diamond}(A) \neq f(A)$ for square indefinite matrices $A$, yet certain relations exist between $f^{\diamond}(A)$ and $f(A)$ \cite{Aurentz19GMF}. Let $q(x)$ be a polynomial, and define $p(x) = q(x^2)x$. We have
\begin{equation}\label{eq:pq_relation}
	q(A^{\dagger}A)A^{\dagger} = q(V \Sigma^2 V^{\dagger}) V \Sigma W^{\dagger} = V q(\Sigma^2) \Sigma W^{\dagger} = V p(\Sigma) W^{\dagger} = p^{\diamond}(A^{\dagger}).
\end{equation}

Recall that when solving the normal equations using CG method, the approximate solution vector $\bm{x}_j$ lies in the Krylov subspace $\bm{x}_0 + \mathcal{K}_j(A^{\dagger}A, A^{\dagger}\bm{b})$, which means that $\bm{x}_j$ can be represented as the generalized matrix polynomial $ q(A^{\dagger}A)A^{\dagger}= p^{\diamond}(A^{\dagger}) $ applied to vector $\bm{b}$ when choosing $\bm{x}_0 = 0$.

\subsection{Chebyshev expansions}
The Chebyshev polynomial $T_j(x)$ of the first kind is a polynomial in $x$ of degree $j$, defined by the relation
\begin{equation}
T_j(x):=\cos(j\eta),
\end{equation}
where $\eta=\arccos(x)$, $x\in[-1,1]$. The Chebyshev polynomials are orthogonal in the sense that
\begin{align*}
\int_{-1}^{1}T_j(x)T_k(x)\frac{dx}{\sqrt{1-x^2}} &= \begin{cases}
0, & j\neq k,  \\
\frac{\pi}{2}, & j=k\neq0,\\
\pi, & j=k=0.
\end{cases}
\end{align*}
Then, a polynomial $f(x)$ of degree $d$ can be expanded in terms of Chebyshev polynomials as
\begin{equation*}
f(x)=\sum\limits_{j=0}^{d}a_jT_j(x),
\end{equation*}
where the expansion coefficients are determined by
\begin{align*}
a_j &= \begin{cases}
\frac{2}{\pi}\int_{-1}^{1}T_j(x)f(x)\frac{dx}{\sqrt{1-x^2}}, & j\geq 1,  \\
\frac{1}{\pi}\int_{-1}^{1}T_0(x)f(x)\frac{dx}{\sqrt{1-x^2}}, & j=0.
\end{cases}
\end{align*}
We apply the substitution $\eta=\arccos(x)$ and re-express the Chebyshev expansion as follows:
\begin{equation}\label{eq:Chebyshev expansion}
f(\cos(\eta))=\sum\limits_{j=0}^{d}a_j\cos(j\eta)
=\sum\limits_{j=0}^{d}a_j\frac{e^{ij\eta}+e^{-ij\eta}}{2}
=\sum\limits_{j=-d}^{d}c_je^{ij\eta}:=\tilde{f}_d(e^{i\eta}),
\end{equation}
where $c_j=a_j/2$ for $j\geq0$ and $c_j=a_{-j}/2$ for $j<0$.
We define $\hat{f}_{2d}(e^{i\eta})=\sum\limits_{j=-d}^{d}c_je^{i(d+j)\eta}$. Then we have $\hat{f}_{2d}(e^{i\eta})=e^{id\eta}\tilde{f}_d(e^{i\eta})$.

\subsection{Generalized quantum signal processing}
Let $H$ be a Hamiltonian. The signal operator used in generalized quantum signal processing is a $0$-controlled application of $U = e^{iH}$:
$$M = \left(
\begin{array}{cc}
	U & 0 \\
	0 & I \\
\end{array}
\right).$$
The signal processing operations are defined as arbitrary $SU(2)$ rotations of the ancillary qubit:
\begin{equation}\label{eq:rot_tensor}
	R(\theta,\phi,\lambda)=\left(
	\begin{array}{cc}
		e^{i(\lambda+\phi)}\cos\theta & e^{i\phi}\sin\theta \\
		e^{i\lambda}\sin\theta & -\cos\theta \\
	\end{array}
	\right)\otimes I.
\end{equation}
For the sake of brevity, we denote 
\begin{equation}\label{eq:rot_def1}
	R_0:=R(\theta_0,\phi_0,\lambda)=\Upsilon_0\otimes I,~~ R_j:=R(\theta_j,\phi_j,0)=\Upsilon_j\otimes I ~(j\geq1),
\end{equation}
where
\begin{equation}\label{eq:rot_def2}
\Upsilon_0:=\left(
\begin{array}{cc}
	e^{i(\lambda+\phi_0)}\cos\theta_0 & e^{i\phi_0}\sin\theta_0 \\
	e^{i\lambda}\sin\theta_0 & -\cos\theta_0 \\
\end{array}
\right)
,~~
\Upsilon_j:=\left(
\begin{array}{cc}
	e^{i\phi_j}\cos\theta_j & e^{i\phi_j}\sin\theta_j \\
	\sin\theta_j & -\cos\theta_j \\
\end{array}
\right).
\end{equation}

The conclusion in Ref.\cite{MW24GQSP} shows that the polynomial transformations of the unitary matrix $U$ can be block encoded by alternately performing the signal operator and signal processing operator. The quantum circuit is shown in Fig.~\ref{fig:GQSP}. Importantly, these polynomials require no additional assumptions after being appropriately scaled, as shown in the following lemma.

\begin{lemma}(Generalized quantum signal processing, GQSP\cite{MW24GQSP})\label{lemma:GQSP}
	$\forall P\in\mathbb{C}[x]$, with $deg(P)=d$, if $\forall x\in\mathbb{R}$, $|P(e^{ix})|^2\leq1$.
	Then $\exists\bm{\theta},\bm{\phi}\in\mathbb{R}^{d+1}$, $\lambda\in\mathbb{R}$ such that
	$$\left(
	\begin{array}{cc}
		P(U) & \cdot \\
		\cdot & \cdot \\
	\end{array}
	\right)=(\prod_{j=1}^{d}R(\theta_j,\phi_j,0)M)R(\theta_0,\phi_0,\lambda).$$
\end{lemma}

\begin{figure}[ht]
	\[
	\Qcircuit @C=.8em @R=1.6em {
		&\lstick{} &\qw    &\gate{\Upsilon_0} &\qw
		&\ctrlo{1} &\qw    &\gate{\Upsilon_1}       &\qw  &\ctrlo{1} &\qw    &\cdots  & & \ctrlo{1}   &\gate{\Upsilon_d} & \qw &\qw
		\\
		&\lstick{} &\qw  &\qw & \qw &\multigate{2}{U}  &\qw & \qw &\qw &\multigate{2}{U} &\qw & \cdots & &\multigate{2}{U}  & \qw &\qw &\qw
		\\
		&\lstick{} &  &\vdots &
		&  & & & &  &  &   & &     &  &   &
		\\
		&\lstick{} &\qw  &\qw  &\qw &\ghost{U} &\qw  &\qw &\qw &\ghost{U}    &\qw & \cdots & & \ghost{U} &\qw & \qw &\qw }
	\]
	\caption{Quantum circuit for GQSP, where the phase factors $\{\theta_i\}_{i=0}^d$, $\{\varphi_i\}_{i=0}^d$ and $\lambda$ are chosen to enact a polynomial transformation of the unitary $U$. Here $\Upsilon_j$ is a general single qubit rotation defined in Eq.(\ref{eq:rot_def2}).}
	\label{fig:GQSP}
\end{figure}

\section{Generalized quantum singular value transformation}\label{sect:GQSVT}
Let the singular value decomposition of matrix $A$ be $A=W_{\Sigma}\Sigma V_{\Sigma}^{\dag}$.
Denote the columns of $W_{\Sigma}$ and $V_{\Sigma}$ by $\{\ket{\bm{w}_k}\}$ and $\{\ket{\bm{v}_k}\}$, respectively. Then we can rewrite $A$ as
$$A = \sum_{k=1}^r \sigma_k \ket{\bm{w}_k} \bra{\bm{v}_k}.$$
Without loss of generality, assume that $\|A\|\leq1$. 

\subsection{Qubitization}\label{sect:Qubitization}
\begin{definition}
	For an $n\times n$ matrix $A$, $n=2^s$, an unitary matrix $U$ is called $(\alpha,a,\varepsilon)$-block encoding of $A$ if there exist $\alpha, \varepsilon\in\mathbb{R}_{+}$ such that
	\begin{equation*}
		\|A-\alpha(\bra{0^a}\otimes I_n)U(\ket{0^a}\otimes I_n))\|\leq\varepsilon.
	\end{equation*}
\end{definition}

Suppose we are given a unitary matrix $U$ that block encodes $A$, i.e.,
\begin{equation}\label{eq:block}
U=\begin{array}{c}
\begin{array} {cc}
              &  \begin{array}{cc}  \Pi & \  \end{array}\\
\begin{array}{c}
\tilde{\Pi}\\
\
\end{array}  & \left(\begin{array}{cc}
		A & \cdot  \\
		\cdot & \cdot \end{array}\right),
\end{array}\\
\
\end{array}
\end{equation}
where $\Pi=\ket{0}\bra{0}\otimes\sum\limits_{k} \ket{\bm{v}_k}\bra{\bm{v}_k}$ and $\tilde{\Pi}=\ket{0}\bra{0}\otimes\sum\limits_{k} \ket{\bm{w}_k}\bra{\bm{w}_k}$ are projectors, such that $\tilde{\Pi}U\Pi=\ket{0}\bra{0}\otimes A$. We can supplement the missing blocks of $U$ as done in Ref.\cite{MRT21Grand} to obtain
\begin{equation*}
U =\left(
      \begin{array}{cc}
        A & \sqrt{I-A^2} \\
        \sqrt{I-A^2} & -A \\
      \end{array}
    \right),
\end{equation*}
where $\sqrt{I-A^2}$ is formally defined as $\sqrt{I-A^2}:=\sum_{k=1}^r \sqrt{1-\sigma_k^2}\ket{\bm{w}_k}\bra{\bm{v}_k}$. Then, we have
\begin{equation*}
U(\ket{0}\ket{\bm{v}_k}, \ket{1}\ket{\bm{v}_k})=(\ket{0}\ket{\bm{w}_k}, \ket{1}\ket{\bm{w}_k})\left(
      \begin{array}{cc}
        \sigma_k & \sqrt{1-\sigma_k^2} \\
        \sqrt{1-\sigma_k^2} & -\sigma_k \\
      \end{array}
    \right),
\end{equation*}
\begin{equation*}
U^{\dag}(\ket{0}\ket{\bm{w}_k}, \ket{1}\ket{\bm{w}_k})=(\ket{0}\ket{\bm{v}_k}, \ket{1}\ket{\bm{v}_k})\left(
      \begin{array}{cc}
        \sigma_k & \sqrt{1-\sigma_k^2} \\
         \sqrt{1-\sigma_k^2} & -\sigma_k \\
      \end{array}
    \right),
\end{equation*}
\begin{equation*}
(2\Pi-I)(\ket{0}\ket{\bm{v}_k}, \ket{1}\ket{\bm{v}_k})=(\ket{0}\ket{\bm{v}_k}, \ket{1}\ket{\bm{v}_k})\left(
      \begin{array}{cc}
        1 & 0 \\
        0 & -1 \\
      \end{array}
    \right),
\end{equation*}
\begin{equation*}
(2\tilde{\Pi}-I)(\ket{0}\ket{\bm{w}_k}, \ket{1}\ket{\bm{w}_k})=(\ket{0}\ket{\bm{w}_k}, \ket{1}\ket{\bm{w}_k})\left(
      \begin{array}{cc}
        1 & 0 \\
        0 & -1 \\
      \end{array}
    \right).
\end{equation*}

Define unitary operators $W=(2\tilde{\Pi}-I)U$ and $\tilde{W}=(2\Pi-I)U^{\dag}$. The action of $W$ restricted to the two-dimensional subspace $\operatorname{span}\{\ket{0}\ket{\bm{v}_k}, \ket{1}\ket{\bm{v}_k}\}$ is
\begin{equation*}
W_{\sigma_k} = \left(
      \begin{array}{cc}
        \sigma_k & \sqrt{1-\sigma_k^2} \\
        -\sqrt{1-\sigma_k^2} & \sigma_k \\
      \end{array}
    \right).
\end{equation*}
Therefore, the operator $W$ can be expressed in the form of a direct sum:
\begin{equation}
W = \sum\limits_{k}\left(
      \begin{array}{cc}
        \sigma_k & \sqrt{1-\sigma_k^2} \\
        -\sqrt{1-\sigma_k^2} & \sigma_k \\
      \end{array}
    \right)\otimes\ket{\bm{w}_k}\bra{\bm{v}_k}=\sum\limits_{k}W_{\sigma_k}\otimes\ket{\bm{w}_k}\bra{\bm{v}_k}.
\end{equation}
Similarly, the action of $\tilde{W}$ restricted to the two-dimensional subspace spanned by $\{\ket{0}\ket{\bm{w}_k}, \ket{1}\ket{\bm{w}_k}\}$ is also $W_{\sigma_k}$, and $\tilde{W}$ has the direct sum form
\begin{equation}
\tilde{W} = \sum\limits_{k}\left(
      \begin{array}{cc}
        \sigma_k & \sqrt{1-\sigma_k^2} \\
        -\sqrt{1-\sigma_k^2} & \sigma_k \\
      \end{array}
    \right)\otimes\ket{\bm{v}_k}\bra{\bm{w}_k}=\sum\limits_{k}W_{\sigma_k}\otimes\ket{\bm{v}_k}\bra{\bm{w}_k}.
\end{equation}

Note that the eigenvectors corresponding to the eigenvalues $e^{\pm i\eta_k}$ of matrix $W_{\sigma_k}$ are $\ket{\bm{\omega}_{\pm}}=\frac{1}{\sqrt{2}}(\ket{0}\pm i\ket{1})$, where $\eta_k=\arccos(\sigma_k)$. Denote $\ket{\bm{\phi}_{k}^{\pm}}=\ket{\bm{\omega}_{\pm}}\ket{\bm{v}_{k}}$,
$\ket{\tilde{\bm{\phi}}_{k}^{\pm}}=\ket{\bm{\omega}_{\pm}}\ket{\bm{w}_{k}}$.
On the two-dimensional subspace spanned by $\{\ket{\bm{\phi}_{k}^{+}}, \ket{\bm{\phi}_{k}^{-}}\}$, the action of $W$ is
\begin{equation*}
W(\ket{\bm{\phi}_{k}^{+}}, \ket{\bm{\phi}_{k}^{-}})=(\ket{\tilde{\bm{\phi}}_{k}^{+}}, \ket{\tilde{\bm{\phi}}_{k}^{-}})\left(
      \begin{array}{cc}
        e^{i\eta_k} & \  \\
        \  & e^{-i\eta_k} \\
      \end{array}
    \right).
\end{equation*}
The action of $\tilde{W}$ on two-dimensional subspace spanned by $\{|\tilde{\bm{\phi}}_{k}^{+}\rangle, |\tilde{\bm{\phi}}_{k}^{-}\rangle\}$ is
\begin{equation*}
\tilde{W}(\ket{\tilde{\bm{\phi}}_{k}^{+}}, \ket{\tilde{\bm{\phi}}_{k}^{-}})=(\ket{\bm{\phi}_{k}^{+}}, \ket{\bm{\phi}_{k}^{-}})\left(
      \begin{array}{cc}
        e^{i\eta_k} & \  \\
        \  & e^{-i\eta_k} \\
      \end{array}
    \right).
\end{equation*}
We see that the operator $W$ brings a state in the space spanned by $\{\ket{\bm{\phi}_{k}^{+}}, \ket{\bm{\phi}_{k}^{-}}\}$ into a state in the space spanned by $\{\ket{\tilde{\bm{\phi}}_{k}^{+}}, \ket{\tilde{\bm{\phi}}_{k}^{-}}\}$, and the operator $\tilde{W}$ does the reverse.

\subsection{Four controlled operators}\label{sect:Four_controlled_operators}
Define $0$-controlled operators 
\begin{equation}\label{eq:M&Mtil}
M=\ket{0}\bra{0}\otimes W+\ket{1}\bra{1}\otimes I ~~\text{and}~~ \tilde{M}=\ket{0}\bra{0}\otimes \tilde{W}+\ket{1}\bra{1}\otimes I.
\end{equation}
We can verify that
\begin{equation*}
M\left(
      \begin{array}{cccc}
        \ket{\bm{\phi}_{k}^{+}}  & \ &\ket{\bm{\phi}_{k}^{-}} &\  \\
        \ &\ket{\bm{\phi}_{k}^{+}} &\  & \ket{\bm{\phi}_{k}^{-}} \\
      \end{array}
    \right)=\left(
      \begin{array}{cccc}
        \ket{\tilde{\bm{\phi}}_{k}^{+}} & \ &\ket{\tilde{\bm{\phi}}_{k}^{-}} &\  \\
        \ &\ket{\tilde{\bm{\phi}}_{k}^{+}} &\  & \ket{\tilde{\bm{\phi}}_{k}^{-}} \\
      \end{array}
    \right)\left(
      \begin{array}{cccc}
        e^{i\eta_k} &\ &\ &\  \\
        \  &1  &\ &\ \\
        \ &\ &e^{-i\eta_k} &\ \\
        \ &\ &\ &1 \\
      \end{array}
    \right).
\end{equation*}
Consequently, the operator $M$ can be expressed as a direct sum
\begin{equation}\label{eq:direct_sum_M}
\begin{split}
M&=\sum\limits_k\left(
      \begin{array}{cccc}
        \ket{\tilde{\bm{\phi}}_{k}^{+}} & \ &\ket{\tilde{\bm{\phi}}_{k}^{-}} &\ \\
        \ &\ket{\tilde{\bm{\phi}}_{k}^{+}} &\  & \ket{\tilde{\bm{\phi}}_{k}^{-}} \\
      \end{array}
    \right)\left(
      \begin{array}{cccc}
        e^{i\eta_k} &\ &\ &\  \\
        \  &1  &\ &\ \\
        \ &\ &e^{-i\eta_k} &\ \\
        \ &\ &\ &1 \\
      \end{array}
    \right)\left(
      \begin{array}{cc}
        \bra{\bm{\phi}_{k}^{+}}  & \ \\
        \ &\bra{\bm{\phi}_{k}^{+}} \\
        \bra{\bm{\phi}_{k}^{-}}  &\ \\
        \ &\bra{\bm{\phi}_{k}^{-}} \\
      \end{array}
    \right)\\
&=\sum\limits_k\left(
      \begin{array}{cc}
        e^{i\eta_k}\ket{\tilde{\bm{\phi}}_{k}^{+}}\bra{\bm{\phi}_{k}^{+}}
        +e^{-i\eta_k}\ket{\tilde{\bm{\phi}}_{k}^{-}}\bra{\bm{\phi}_{k}^{-}} & \  \\
        \  & \ket{\tilde{\bm{\phi}}_{k}^{+}}\bra{\bm{\phi}_{k}^{+}}
        +\ket{\tilde{\bm{\phi}}_{k}^{-}}\bra{\bm{\phi}_{k}^{-}} \\
      \end{array}
    \right)\\
&=\sum\limits_k\left(
      \begin{array}{cc}
        e^{i\eta_k} & \  \\
        \  & 1 \\
      \end{array}
    \right)\otimes\ket{\tilde{\bm{\phi}}_{k}^{+}}\bra{\bm{\phi}_{k}^{+}}
    +\left(
      \begin{array}{cc}
        e^{-i\eta_k} & \  \\
        \  & 1 \\
      \end{array}
    \right)\otimes\ket{\tilde{\bm{\phi}}_{k}^{-}}\bra{\bm{\phi}_{k}^{-}}.
\end{split}
\end{equation}
In a similar manner, one can express $\tilde{M}$ as a direct sum
\begin{equation}\label{eq:direct_sum_Mtil}
\tilde{M}=\sum\limits_k\left[\left(
      \begin{array}{cc}
        e^{i\eta_k} & \  \\
        \  & 1 \\
      \end{array}
    \right)\otimes\ket{\bm{\phi}_{k}^{+}}\bra{\tilde{\bm{\phi}}_{k}^{+}}
    +\left(
      \begin{array}{cc}
        e^{-i\eta_k} & \  \\
        \  & 1 \\
      \end{array}
    \right)\otimes\ket{\bm{\phi}_{k}^{-}}\bra{\tilde{\bm{\phi}}_{k}^{-}}\right].
\end{equation}

We define $1$-controlled operators
\begin{equation}\label{eq:N&Ntil}
N=|0\rangle\langle0|\otimes I+|1\rangle\langle1|\otimes W^{\dagger}~~\text{and}~~ \tilde{N}=|0\rangle\langle0|\otimes I +|1\rangle\langle1|\otimes \tilde{W}^{\dagger}.
\end{equation}
Similarly with the process as described above, $N$ and $\tilde{N}$ are presented in direct sum forms
\begin{equation}\label{eq:direct_sum_N}
N=\sum\limits_k\left[\left(
      \begin{array}{cc}
        1  & \  \\
        \  &  e^{-i\eta_k} \\
      \end{array}
    \right)\otimes\ket{\bm{\phi}_{k}^{+}}\bra{\tilde{\bm{\phi}}_{k}^{+}}
    +\left(
      \begin{array}{cc}
        1  & \  \\
        \  & e^{i\eta_k}\\
      \end{array}
    \right)\otimes\ket{\bm{\phi}_{k}^{-}}\bra{\tilde{\bm{\phi}}_{k}^{-}}\right],
\end{equation}
\begin{equation}\label{eq:direct_sum_Ntil}
\tilde{N}=\sum\limits_k\left[\left(
      \begin{array}{cc}
        1  & \  \\
        \  &  e^{-i\eta_k} \\
      \end{array}
    \right)\otimes\ket{\tilde{\bm{\phi}}_{k}^{+}}\bra{\bm{\phi}_{k}^{+}}
    +\left(
      \begin{array}{cc}
        1  & \  \\
        \  & e^{i\eta_k}\\
      \end{array}
    \right)\otimes\ket{\tilde{\bm{\phi}}_{k}^{-}}\bra{\bm{\phi}_{k}^{-}}\right].
\end{equation}

Next, we will use the frame of GQSP to perform generalized signal processing for general matrices. The signal processing operations we use are consistent with GQSP, as shown in Eq.(\ref{eq:rot_tensor}). The signal operators we use are $M$, $\tilde{M}$, $N$ and $\tilde{N}$, as defined in Eq.(\ref{eq:M&Mtil}) and Eq.(\ref{eq:N&Ntil}). Our method requires one more ancilla qubit than GQSP, namely two ancilla qubits are needed.

\subsection{Polynomial transform of singular values}
For a target polynomial $f(x)$ of degree $d$, we have the corresponding Chebyshev expansion $\tilde{f}_d(e^{i\eta})$ as described in Eq.(\ref{eq:Chebyshev expansion}). Then, by using Lemma \ref{lemma:GQSP}, there exist $\bm{\theta},\bm{\phi}\in \mathbb{R}^{2d+1}$ and $\lambda\in \mathbb{R}$ for $\hat{f}_{2d}(e^{i\eta})=e^{id\eta}\tilde{f}_d(e^{i\eta})$ such that
\begin{equation}\label{eq:subspace}
\prod\limits_{j=1}^{2d}\left[\Upsilon_j\left(
      \begin{array}{cc}
        e^{i\eta} & \  \\
        \  & 1 \\
      \end{array}
    \right)\right]\Upsilon_0=\left(\begin{array}{cc}
        \hat{f}_{2d}(e^{i\eta}) & \cdot \\
        \cdot & \cdot \\
      \end{array}\right).
\end{equation}

We first discuss the polynomial transformations on individual singular value. Assume that the initial state is given by
$$\ket{0}\ket{0}\ket{\bm{v}_k},$$
where $\ket{\bm{v}_k}$ is the singular vector corresponding to the singular value $\sigma_k$. Let $\eta_k=\arccos(\sigma_k)$. 

We take the case where $d$ is even as an example. By successively applying the operators $R_0$, $M$, $R_1$, $\tilde{M}$, $R_2$, $\ldots$, $\tilde{M}$ and $R_d$, together with the direct sum representations of $M$ and $\tilde{M}$ in Eqs.~(\ref{eq:direct_sum_M})-(\ref{eq:direct_sum_Mtil}), we obtain
\begin{equation}\label{eq:operators_R_M}
\begin{split}
\ket{0}\ket{0}\ket{\bm{v}_k}~=&~
\frac{1}{\sqrt{2}}\ket{0}\ket{\bm{\phi}_k^{+}}+\frac{1}{\sqrt{2}}\ket{0}\ket{\bm{\phi}_k^{-}}\\
\xrightarrow{R_0}&~
\frac{1}{\sqrt{2}}\Upsilon_0 \ket{0}\ket{\bm{\phi}_k^{+}}+\frac{1}{\sqrt{2}}\Upsilon_0 \ket{0}\ket{\bm{\phi}_k^{-}}\\
\xrightarrow{M,R_1} &~
\frac{1}{\sqrt{2}}\Upsilon_1\left(
      \begin{array}{cc}
        e^{i\eta_k} & \  \\
        \  & 1 \\
      \end{array}
    \right)\Upsilon_0 \ket{0}\ket{\tilde{\bm{\phi}}_k^{+}}
+\frac{1}{\sqrt{2}}\Upsilon_1\left(
      \begin{array}{cc}
        e^{-i\eta_k} & \  \\
        \  & 1 \\
      \end{array}
    \right)\Upsilon_0 \ket{0}\ket{\tilde{\bm{\phi}}_k^{-}}\\
\xrightarrow{\tilde{M},R_2} &~
\frac{1}{\sqrt{2}}\Upsilon_2\left(
      \begin{array}{cc}
        e^{i\eta_k} & \  \\
        \  & 1 \\
      \end{array}
    \right)\Upsilon_1\left(
      \begin{array}{cc}
        e^{i\eta_k} & \  \\
        \  & 1 \\
      \end{array}
    \right)\Upsilon_0 \ket{0}\ket{\bm{\phi}_k^{+}}\\
&+\frac{1}{\sqrt{2}}\Upsilon_2\left(
      \begin{array}{cc}
        e^{-i\eta_k} & \  \\
        \  & 1 \\
      \end{array}
    \right)\Upsilon_1\left(
      \begin{array}{cc}
        e^{-i\eta_k} & \  \\
        \  & 1 \\
      \end{array}
    \right)\Upsilon_0 \ket{0}\ket{\bm{\phi}_k^{-}}\\
\longrightarrow &~ \cdots\\
\xrightarrow{\tilde{M},R_d}&~\frac{1}{\sqrt{2}}\prod\limits_{j=1}^{d}\left[\Upsilon_j\left(
      \begin{array}{cc}
        e^{i\eta_k} & \  \\
        \  & 1 \\
      \end{array}
    \right)\right]\Upsilon_0 \ket{0}\ket{\bm{\phi}_k^{+}}
+\frac{1}{\sqrt{2}}\prod\limits_{j=1}^{d}\left[\Upsilon_j\left(
      \begin{array}{cc}
        e^{-i\eta_k} & \  \\
        \  & 1 \\
      \end{array}
    \right)\right]\Upsilon_0\ket{0}\ket{\bm{\phi}_k^{-}}.
\end{split}
\end{equation}

As observed in Eq.(\ref{eq:Chebyshev expansion}), the polynomial $f(x)$ is expressed as a Laurent polynomial in $e^{i\eta}$, incorporating terms with negative powers. Consequently, we apply the operators $\tilde{N}$, $R_{d+1}$, $N$, $R_{d+2}$, $\ldots$, $N$, and $R_{2d}$ in sequence to introduce negative power terms, leading to the result
\begin{equation}\label{eq:operators_R_N}
\begin{split}
\xrightarrow{\tilde{N},R_{d+1}}&~  \frac{1}{\sqrt{2}}\Upsilon_{d+1}
    \left(\begin{array}{cc}
        1 & \  \\
        \  & e^{-i\eta_k} \\
      \end{array}
    \right)\prod\limits_{j=1}^{d}\left[\Upsilon_j\left(
      \begin{array}{cc}
        e^{i\eta_k} & \  \\
        \  & 1 \\
      \end{array}
    \right)\right]\Upsilon_0 \ket{0}\ket{\tilde{\bm{\phi}}_k^{+}}\\
&+\frac{1}{\sqrt{2}}\Upsilon_{d+1}
\left(\begin{array}{cc}
        1 & \  \\
        \  & e^{i\eta_k} \\
      \end{array}
    \right)\prod\limits_{j=1}^{d}\left[\Upsilon_j\left(
      \begin{array}{cc}
        e^{-i\eta_k} & \  \\
        \  & 1 \\
      \end{array}
    \right)\right]\Upsilon_0 \ket{0}\ket{\tilde{\bm{\phi}}_k^{-}}\\
\longrightarrow &~ \cdots\\
\xrightarrow{N,R_{2d}}&~ \frac{1}{\sqrt{2}}\Upsilon_{2d}
\left(\begin{array}{cc}
        1 & \  \\
        \  & e^{-i\eta_k}  \\
      \end{array}
    \right)\cdots \Upsilon_{d+1}
\left(\begin{array}{cc}
        1 & \  \\
        \  & e^{-i\eta_k} \\
      \end{array}
    \right)
    \prod\limits_{j=1}^{d}\left[\Upsilon_j\left(
      \begin{array}{cc}
        e^{i\eta_k} & \  \\
        \  & 1 \\
      \end{array}
    \right)\right]\Upsilon_0 \ket{0}\ket{\bm{\phi}_k^{+}}\\
&+\frac{1}{\sqrt{2}}\Upsilon_{2d}
\left(\begin{array}{cc}
        1 & \  \\
        \  & e^{i\eta_k} \\
      \end{array}
    \right)
    \cdots \Upsilon_{d+1}
\left(\begin{array}{cc}
        1 & \  \\
        \  & e^{i\eta_k} \\
      \end{array}
    \right)\prod\limits_{j=1}^{d}\left[\Upsilon_j\left(
      \begin{array}{cc}
        e^{-i\eta_k} & \  \\
        \  & 1 \\
      \end{array}
    \right)\right]\Upsilon_0 \ket{0}\ket{\bm{\phi}_k^{-}}\\
=& \frac{1}{\sqrt{2}}e^{-id\eta_k}\prod\limits_{j=1}^{2d}\left[\Upsilon_j\left(
      \begin{array}{cc}
        e^{i\eta_k} & \  \\
        \  & 1 \\
      \end{array}
    \right)\right]\Upsilon_0 \ket{0}\ket{\bm{\phi}_k^{+}}+\frac{1}{\sqrt{2}}
e^{id\eta_k}\prod\limits_{j=1}^{2d}\left[\Upsilon_j\left(
      \begin{array}{cc}
        e^{-i\eta_k} & \  \\
        \  & 1 \\
      \end{array}
    \right)\right]\Upsilon_0 \ket{0}\ket{\bm{\phi}_k^{-}}.
\end{split}
\end{equation}

Using Eq.~(\ref{eq:subspace}), the above expression can be further simplified to
\begin{equation}
\begin{split}
&\frac{1}{\sqrt{2}}e^{-id\eta_k}\left(\begin{array}{cc}
        \hat{f}_{2d}(e^{i\eta_k}) & \cdot \\
        \cdot & \cdot \\
      \end{array}
    \right)\ket{0}\ket{\bm{\phi}_k^{+}}
+\frac{1}{\sqrt{2}}e^{id\eta_k}\left(\begin{array}{cc}
      \hat{f}_{2d}(e^{-i\eta_k}) & \cdot \\
        \cdot & \cdot\\
      \end{array}
    \right)\ket{0}\ket{\bm{\phi}_k^{-}}\\
=&~\ket{0}\left[\frac{1}{\sqrt{2}}e^{-id\eta_k}\hat{f}_{2d}(e^{i\eta_k})\ket{\bm{\phi}_k^{+}}
+\frac{1}{\sqrt{2}}e^{id\eta_k}\hat{f}_{2d}(e^{-i\eta_k})\ket{\bm{\phi}_k^{-}}\right]+\ket{\perp}\\
=&~\ket{0} \left[\frac{1}{\sqrt{2}}\tilde{f}_d(e^{i\eta_k})\ket{\bm{\phi}_k^{+}}
+\frac{1}{\sqrt{2}}\tilde{f}_d(e^{-i\eta_k})\ket{\bm{\phi}_k^{-}}\right]+\ket{\perp}\\
=&~\ket{0} f(\sigma_k)\ket{0}\ket{\bm{v}_k}+\ket{\perp},\\
\end{split}
\end{equation}
where $\ket{\perp}$ is a state orthogonal to the previous item. 

The above discussion explains how our method implements a polynomial transformation of a single singular value, i.e.
\begin{equation*}
    \ket{0}\ket{0}\ket{\bm{v}_k}\longrightarrow f(\sigma_k)\ket{0} \ket{0}\ket{\bm{v}_k}+\ket{\perp}.
\end{equation*}
The polynomial transformation of the entire matrix will be presented in the next subsection.

\subsection{Matrix polynomial transform}
We assume that the input vector is $\ket{0}\ket{0}\ket{\bm{\varphi}}$, where $\ket{\bm{\varphi}}=\sum_k\varphi_k\ket{\bm{v}_k}$. For the case where the polynomial  $f(x)$ has even degree $d$, applying the same operators as in Eq.~(\ref{eq:operators_R_M}) and Eq.~(\ref{eq:operators_R_N}), we obtain that
\begin{equation}
\ket{0}\ket{0}\ket{\bm{\varphi}}
\longrightarrow\ket{0}\ket{0}\sum_k f(\sigma_k)\varphi_k\ket{\bm{v}_k}+\ket{\perp'}=\ket{0}\ket{0} f^{R}(A)\ket{\bm{\varphi}}+\ket{\perp'},
\end{equation}
where $\ket{\perp'}$ is a state orthogonal to previous item, $f^{R}(A)$ is right generalized matrix polynomial defined as $f^{R}(A)\equiv\sum\limits_kf(\sigma_k)\ket{\bm{v}_k}\bra{\bm{v}_k}$. 
By measuring the first two ancilla qubits, we obtain $f^{R}(A)\ket{\bm{\varphi}}$ when the result is $\ket{00}$. Similarly, we can obtain $f^{L}(A)\ket{\bm{\varphi}}$ and $f^{\diamond}(A)\ket{\bm{\varphi}}$. Here, the left generalized matrix polynomial $f^{L}(A)$ and the generalized matrix polynomial $f^{\diamond}(A)$ are defined  as $f^{L}(A) = \sum_kf(\sigma_k)\ket{\bm{w}_k}\bra{\bm{w}_k}$ and  $f^{\diamond}(A)\equiv\sum\limits_kf(\sigma_k)\ket{\bm{w}_k}\bra{\bm{v}_k}$, respectively. 

We summarize the above discussion in the following theorem.

\begin{theorem}(Generalized Quantum Singular Value Transformation, GQSVT)
Given a block encoding of a matrix $A=\sum_{k=1}^r \sigma_k\ket{\bm{w}_k}\bra{\bm{v}_k}$ in a unitary matrix $U$ as shown in (\ref{eq:block}). The definition of operations $R_0$, $R_j(j\geq1)$, $M$, $\tilde{M}$, $N$ and $\tilde{N}$ are given in Eqs.~(\ref{eq:rot_tensor}), (\ref{eq:M&Mtil}) and (\ref{eq:N&Ntil}), respectively. For any polynomial $f(x)$ of degree $d$, if $\forall x\in[-1,1]$, $|f(x)|\leq1$, then $\exists \bm{\theta}, \bm{\phi}\in \mathbb{R}^{d+1}$, $\lambda\in
\mathbb{R}$, such that
 for even $d$,
\begin{equation}\label{eq:U_even_R}
U_{\bm{\theta},\bm{\phi},\lambda}=
\left(\prod\limits_{j=1}^{d/2}R_{d+2j}N R_{d+2j-1}\tilde{N}\right)\left(\prod\limits_{j=1}^{d/2}R_{2j}\tilde{M}R_{2j-1}M\right)R_0
\end{equation}
is a block encoding of $f^{R}(A)$. 
For odd $d$,
\begin{equation}\label{eq:U_diamond}
U_{\bm{\theta},\bm{\phi},\lambda}=
R_{2d}N\left(\prod\limits_{j=1}^{(d-1)/2}R_{d+2j}\tilde{N}R_{d+2j-1}N\right)
R_{d}M\left(\prod\limits_{j=1}^{(d-1)/2}R_{2j}\tilde{M}R_{2j-1}M\right)R_0
\end{equation}
is a block encoding of $f^{\diamond}(A)$. 
Namely,
\begin{align}
\left(\bra{0}\bra{0}\otimes I\right)U_{\bm{\theta},\bm{\phi},\lambda}\left(\ket{0}\ket{0}\otimes I\right)
=\begin{cases}
\sum\limits_kf(\sigma_k)\ket{\bm{v}_k}\bra{\bm{v}_k}=f^R(A), & d~ \textmd{even},\\
\sum\limits_kf(\sigma_k)\ket{\bm{w}_k}\bra{\bm{v}_k}=f^{\diamond}(A), & d~\textmd{odd}.
\end{cases}
\end{align}
\end{theorem}

{\bfseries Remark 1.} The SVD of the conjugate transpose of $A$ is given by
$A^{\dagger}=\sum_{k=1}^r \sigma_k\ket{\bm{v}_k}\bra{\bm{w}_k}$. The left generalized matrix polynomial 
$$f^{L}(A) = \sum_kf(\sigma_k)\ket{\bm{w}_k}\bra{\bm{w}_k} = f^{R}(A^{\dagger}).$$
Therefore, the left generalized matrix polynomial transform of matrix $A$ is included in the discussion of the GQSVT of $A^{\dagger}$, which is omitted in the theorem above.

To implement the polynomial transform of $A^{\dagger}$, it is sufficient to modify the operator order, as stated in the following corollary.

\begin{corollary} 
For any polynomial $f(x)$ of degree $d$, if $\forall x\in[-1,1]$, $|f(x)|\leq1$, then $\exists \bm{\theta}, \bm{\phi}\in \mathbb{R}^{d+1}$, $\lambda\in
\mathbb{R}$, such that for even $d$,
$$U_{\bm{\theta},\bm{\phi},\lambda}=\left(\prod\limits_{j=1}^{d/2}R_{d+2j}\tilde{N}R_{d+2j-1}N\right)\left(\prod\limits_{j=1}^{d/2}R_{2j}MR_{2j-1}\tilde{M}\right)R_0$$
is a block encoding of $f^{R}(A^{\dagger})$. For odd $d$,
$$U_{\bm{\theta},\bm{\phi},\lambda} = R_{2d}\tilde{N}\left(\prod\limits_{j=1}^{(d-1)/2}R_{d+2j}NR_{d+2j-1}\tilde{N}\right)R_{d}\tilde{M}\left(\prod\limits_{j=1}^{(d-1)/2}R_{2j}MR_{2j-1}\tilde{M}\right)R_0$$
is a block encoding of $f^{\diamond}(A^{\dagger})$.
\end{corollary}


\subsection{Circuit implementation}
This section discusses the specific implementation of the operator sequence, with subscripts $a$ and $s$ added for clarity.

Suppose $U$ is the $(1,a,0)$-block encoding of matrix $A$, and the dimension of $A$ is $n=2^s$. From the definition in Section \ref{sect:Qubitization} and \ref{sect:Four_controlled_operators}, we have
$$U=\ket{0^a}\bra{0^a}\otimes A+\cdots,~~\Pi=\ket{0^a}\bra{0^a}\otimes\sum\limits_{k}\ket{\bm{v}_k}\bra{\bm{v}_k},$$
$$\tilde{\Pi}=\ket{0^a}\bra{0^a}\otimes\sum\limits_{k}\ket{\bm{w}_k}\bra{\bm{w}_k},~~W=(2\tilde{\Pi}-I)U,~~\tilde{W}=(2\Pi-I)U^{\dagger}.$$
Here, $I$ is the $2^{a+s}$ dimensional identity matrix. In this section, all identity matrices are $2^{a+s}$ dimensional unless otherwise specified. 

The quantum circuit corresponding to the operator sequence shown in Eq.~(\ref{eq:U_even_R}) is given in Fig.~\ref{fig:GQSVT}, where
$$M=\ket{0}\bra{0}\otimes W+\ket{1}\bra{1}\otimes I,~~
\tilde{M}=\ket{0}\bra{0}\otimes \tilde{W}+\ket{1}\bra{1}\otimes I,$$
$$N=|0\rangle\langle0|\otimes I+|1\rangle\langle1|\otimes W^{\dag},~~ \tilde{N}=|0\rangle\langle0|\otimes I +|1\rangle\langle1|\otimes \tilde{W}^{\dag}.$$
\begin{figure}[ht]
\[
\Qcircuit @C=1.0em @R=1.6em {
&\lstick{}&\qw &\gate{\Upsilon_0} &\multigate{2}{M}  &\gate{\Upsilon_1}& \qw & \cdots &  &\multigate{2}{\tilde{M}} &\gate{\Upsilon_d} &\multigate{2}{\tilde{N}} &\gate{\Upsilon_{d+1}}& \qw &\cdots & & \multigate{2}{N}   &\gate{\Upsilon_{2d}} & \qw
\\
&\lstick{} &\qw{/^a} &\qw &\ghost{M}  &\qw & \qw& \cdots & &\ghost{\tilde{M}}  &\qw &\ghost{\tilde{N}} & \qw & \qw & \cdots & &\ghost{N}  & \qw &\qw
\\
&\lstick{}   &\qw{/^s} &\qw &\ghost{M}  &\qw & \qw & \cdots & &\ghost{\tilde{M}}  &\qw &\ghost{\tilde{N}} & \qw & \qw & \cdots & &\ghost{N}  & \qw &\qw
}
\]
\caption{The GQSVT circuit for block encoding a right generalized matrix polynomial $f^{R}(A)$ of degree $d$, where $d$ is even.}
\label{fig:GQSVT} 
\end{figure}

Define the projector-controlled-NOT gate
$$C_{\Pi}\text{NOT}=X\otimes\Pi+I_2\otimes(I-\Pi).$$
The operator $M$ can be implemented by using controlled-$U$, Pauli $Z$, and $C_{\Pi}\text{NOT}$, as shown in Fig.~\ref{fig:3}. Denote the operators in the dashed box in Fig.\ref{fig:M} by $\tilde{\Pi}_Z=C_{I-\tilde{\Pi}}\text{NOT}(Z\otimes I)C_{I-\tilde{\Pi}}\text{NOT}$. We have
\begin{equation*}
	\begin{split}
		\tilde{\Pi}_Z=&C_{I-\tilde{\Pi}}\text{NOT}(Z\otimes I)C_{I-\tilde{\Pi}}\text{NOT}\\
		=&Z\otimes\tilde{\Pi}+(XZX)\otimes(I-\tilde{\Pi})\\
		=&\left(
		\begin{array}{cc}
			2\tilde{\Pi}-I & \  \\
			\  & I-2\tilde{\Pi} \\
		\end{array}
		\right)\\
		=&\ket{0}\bra{0}\otimes(2\tilde{\Pi}-I)+\ket{1}\bra{1}\otimes(I- 2\tilde{\Pi}).
	\end{split}
\end{equation*}
When the first register is $\ket{0}$, $\tilde{\Pi}_Z$ is equivalent to  $\ket{0}\bra{0}\otimes(2\tilde{\Pi}-I)+\ket{1}\bra{1}\otimes I$. Therefore, the quantum circuit on the right side in Fig.\ref{fig:M} is
\begin{equation*}
	\begin{split}
		&\left[\ket{0}\bra{0}\otimes(2\tilde{\Pi}-I)+\ket{1}\bra{1}\otimes I\right]\cdot(\ket{0}\bra{0}\otimes U+\ket{1}\bra{1}\otimes I)\\
		=&\ket{0}\bra{0}\otimes(2\tilde{\Pi}-I)U+\ket{1}\bra{1}\otimes I\\
		=&\ket{0}\bra{0}\otimes W+\ket{1}\bra{1}\otimes I,
	\end{split}
\end{equation*}
which is equivalent to the left side in Fig.\ref{fig:M}. The circuit implementation for $\tilde{M}$, $N$, and $\tilde{N}$ is similar and is therefore omitted here.

\begin{figure}[h!]
\centering
\begin{subfigure}[b]{0.5\textwidth}
\[
\Qcircuit @C=0.5em @R=0.4em @!R {
&\lstick{} & \qw & \qw& \multigate{3}{M} &\qw & \qw & & &\qw & \qw &\ctrlo{1} &\targ & \gate{Z} & \targ &\qw  \\
&\lstick{} &\qw{/^a}& \qw  &\ghost{M} &\qw & \qw  & & &\qw{/^a}& \qw &\multigate{2}{U} & \multigate{2}{I-\tilde{\Pi}}\qwx &\qw &\multigate{2}{I-\tilde{\Pi}}\qwx &\qw \\
&   &  &  &\nghost{M} &  & &\push{\rule{.1em}{0em}=\rule{.1em}{0em}} & & & & &\nghost{U} &\nghost{I-\tilde{\Pi}}   & & \nghost{I-\tilde{\Pi}} & \\
&\lstick{}  &\qw{/^s} & \qw &\ghost{M} & \qw &\qw & & &\qw{/^s}& \qw &\ghost{U} &\ghost{I-\tilde{\Pi}}  &\qw &\ghost{I-\tilde{\Pi}} &\qw
\gategroup{1}{13}{4}{15}{.5em, 1em}{--}
}
\]
\caption{}
\label{fig:M}
\end{subfigure}
\vfill
\begin{subfigure}[b]{0.5\textwidth}
\[
\Qcircuit @C=0.5em @R=0.8em @!R {
&\lstick{} &\qw  & \qw &\targ &\qw & \qw& & & \qw&\qw &\targ &\targ&\targ&\targ&\targ&\targ &\targ &\qw  \\
&\lstick{} &\qw & \qw &\multigate{3}{I-\tilde{\Pi}}\qwx &\qw& \qw & & &\qw & \qw&\ctrl{-1} &\ctrlo{-1}&\ctrlo{-1}&\ctrl{-1}&\ctrlo{-1}&\ctrl{-1}&\ctrl{-1}&\qw \\
&\lstick{} &\qw & \qw &\ghost{I-\tilde{\Pi}}  &\qw & \qw &\push{\rule{.3em}{0em}=\rule{.3em}{0em}} &  &\qw& \qw & \ctrlo{-1} &\ctrl{-1}&\ctrlo{-1}&\ctrl{-1}&\ctrl{-1}&\ctrlo{-1}&\ctrl{-1} &\qw\\
&\lstick{}  &\qw & \qw&\ghost{I-\tilde{\Pi}} & \qw&\qw & & &\qw& \qw &\ctrlo{-1} &\ctrlo{-1}&\ctrl{-1}&\ctrlo{-1}&\ctrl{-1}&\ctrl{-1}&\ctrl{-1} & \qw \\
&\lstick{}  &\qw{/^s}  & \qw&\ghost{I-\tilde{\Pi}} &\qw & \qw& & & \qw{/^s} & \qw&\qw &\qw&\qw&\qw&\qw&\qw&\qw & \qw
}
\]
\caption{}
\label{fig:CNOT}
\end{subfigure}
\caption{(a) the detailed circuit of operator $M$, (b) an example for $C_{I-\tilde{\Pi}}\text{NOT}$ when $a=3$. }
\label{fig:3}
\end{figure}

\section{Solving LS problem using the generalized quantum singular value transformation}
Consider a linear system of equations having $n$ unknowns but $m\geq n $ equations. Symbolically, for matrix $A\in \mathbb{R}^{m\times n}$ and vector $\bm{b}\in \mathbb{R}^{m}$, the LS problem \cite{Trefethen97Numerical} is to find a vector $\bm{x}_{*}\in \mathbb{R}^{n}$ such that
\begin{equation*}
\bm{x}_{*}= \arg\min\limits_{\bm{x}\in \mathbb{R}^{n}} \|A\bm{x}-\bm{b}\|.
\end{equation*}
In the case where $m=n$ and $A$ is nonsingular, the solution is uniquely given by $\bm{x}_{*} = A^{-1}\bm{b}$. When $m>n$, meaning there are more equations than unknowns, the system is $\textit{overdetermined}$. In such cases, the least squares method finds a best-fit solution.

\subsection{Conjugate gradient for least squares}
One classical way to solve LS problems is to solve the normal equations
\begin{equation}\label{eq:normal_equation}
A^{\top}A\bm{x}=A^{\top}\bm{b}.
\end{equation}
If $A$ has full rank, this yields a square and positive definite system. In such cases, the CG method can be applied to solve it, i.e., the CGLS method. Algorithm \ref{algorithm.1} presents the steps of the CGLS method.

\begin{algorithm}[ht]
	\caption{CGLS Algorithm.}
	\label{algorithm.1}
	\textbf{Procedure:} 
	\begin{algorithmic}[1]
		\STATE Choose initial guess $\bm{x}_0$, $\bm{r}_0=A^{\top}\bm{b}-A^{\top}A\bm{x}_0$, $\bm{p}_0=\bm{r}_0$
		\FOR{$j=0,1,...$}
		\STATE $\alpha_j=\frac{\langle \bm{r}_j,\bm{r}_j\rangle}{\langle \bm{p}_j,A^{\top}A\bm{p}_j\rangle}$
		\STATE $\bm{x}_{j+1}=\bm{x}_j+\alpha_j\bm{p}_j$
		\STATE $\bm{r}_{j+1}=\bm{r}_j-\alpha_jA^{\top}A\bm{p}_j$
		\IF{$\|\bm{r}_{j+1}\|\leq \epsilon$}
		\STATE  STOP
		\ENDIF 
		\STATE $\beta_j=\frac{\langle \bm{r}_{j+1},\bm{r}_{j+1}\rangle}{\langle \bm{r}_j,\bm{r}_j\rangle}$
		\STATE $\bm{p}_{j+1}=\bm{r}_{j+1}+\beta_j\bm{p}_j$
		\ENDFOR
	\end{algorithmic}	
\end{algorithm}

Given the initial guess
$\bm{x}_0=0$ in Algorithm \ref{algorithm.1}, the $j$-th iterate 
$\bm{x}_j$ belongs to the Krylov subspace $\mathcal{K}_j(A^{\top}A, A^{\top}\bm{b})$. Moreover, the residuals and search directions satisfy $$\operatorname{span}\{\bm{r}_0,\dots,\bm{r}_{j}\}=\operatorname{span}\{\bm{p}_0,\dots,\bm{p}_{j}\}=\mathcal{K}_{j+1}(A^{\top}A,A^{\top}\bm{b}).$$
Therefore, the vectors produced by Algorithm \ref{algorithm.1} can be expressed as
\begin{equation}\label{eq:poly}
	\begin{split}
		\bm{x}_j=\mathcal{X}_j(A^{\top}A)A^{\top}\bm{b},~~~
		\bm{r}_j=\mathcal{R}_j(A^{\top}A)A^{\top}\bm{b},~~~
		\bm{p}_j=\mathcal{P}_j(A^{\top}A)A^{\top}\bm{b},
	\end{split}
\end{equation}
where $\mathcal{X}_j$, $\mathcal{R}_j$, and $\mathcal{P}_j$ are polynomials defined by 
Algorithm \ref{algorithm.1}. We can verify that these polynomials satisfy the following recursion relations:
$$\mathcal{R}_{j+1}(t)=\mathcal{R}_j(t)-\alpha_jt\mathcal{P}_j(t),~~\mathcal{P}_{j+1}(t)=\mathcal{R}_{j+1}(t)+\beta_j\mathcal{P}_j(t).$$

Defined  $\mathcal{X}_j(t)=\sum_{l=0}^{j-1}\chi_l^{(j)}t^l$, $\mathcal{R}_j(t)=\sum_{l=0}^{j}\gamma_l^{(j)}t^l$ and  $\mathcal{P}_j(t)=\sum_{l=0}^{j}\rho_l^{(j)}t^l$.
Then, the coefficients of these polynomials can be obtained as follows. For $l=0,1,\cdots,j+1$,
\begin{equation*}
\begin{split}
&\chi_{l}^{(j+1)}=\chi_l^{(j)}+\alpha_j\rho_l^{(j)},\\
&\gamma_{l}^{(j+1)}=\gamma_l^{(j)}-\alpha_j\rho_{l-1}^{(j)},\\
&\rho_{l}^{(j+1)}=\gamma_{l}^{(j+1)}+\beta_j\rho_l^{(j)},
\end{split}
\end{equation*}
where $\chi_j^{(j)}=\chi_{j+1}^{(j)}=\gamma_{j+1}^{(j)}=\rho_{-1}^{(j)}=\rho_{j+1}^{(j)}=0$.

According to the relation between standard matrix function and generalized matrix function given in Eq.(\ref{eq:pq_relation}), the iterative vectors in Eq.(\ref{eq:poly}) can be reformulated as the action of the generalized matrix function of $A^{\top}$ on the vector $\bm{b}$. Define
$$\tilde{\mathcal{X}}_j(t) = \mathcal{X}_j(t^2)t,~~ \tilde{\mathcal{R}}_j(t) = \mathcal{R}_j(t^2)t,~~ \tilde{\mathcal{P}}_j(t) = \mathcal{P}_j(t^2)t. $$ 
Then the vectors produced by Algorithm \ref{algorithm.1} be expressed as
\begin{equation}\label{eq:gmp_vectors}
\begin{split}
\bm{x_j} = \tilde{\mathcal{X}}^{\diamond}_j(A^{\top})\bm{b} = \sum_{l=0}^{j-1}\chi_l^{(j)}(A^{\top})^{2l+1}\bm{b},\\
\bm{r_j} = \tilde{\mathcal{R}}^{\diamond}_j(A^{\top})\bm{b} = \sum_{l=0}^{j}\gamma_l^{(j)}(A^{\top})^{2l+1}\bm{b},\\
\bm{p_j} = \tilde{\mathcal{P}}^{\diamond}_j(A^{\top})\bm{b} = \sum_{l=0}^{j}\rho_l^{(j)}(A^{\top})^{2l+1}\bm{b}.
\end{split}
\end{equation}
These representations enable us to incorporate GQSVT into the CGLS algorithm.

\subsection{Hadamard test}\label{sect:swap-test}
In steps $3$ and $9$ of the CGLS algorithm, we need to calculate the inner product. Since we consider a real system, it suffices to compute the real part of the inner product. In this subsection, we introduce the Hadamard test technique for estimating the inner product, which has been employed in a previous study \cite{TWY2024QCG}.

Consider the operators $A$ and $B$ such that $A\ket{0^s}=\ket{\bm{\varphi}}$ and $B\ket{0^s}=\ket{\bm{\psi}}$. Let $U$ and $V$ be $(1,a,0)$ block encodings of $A$ and $B$, respectively, acting on the same qubits as follows:
\begin{equation*}
U\ket{0^a}\ket{0^s}=\ket{0^a}\ket{\bm{\varphi}}+\ket{\perp},~~
V\ket{0^a}\ket{0^s}=\ket{0^a}\ket{\bm{\psi}}+\ket{\perp'},
\end{equation*}
where $\ket{\perp}$ and $\ket{\perp'}$ are states orthogonal to $\ket{0^a}$. The Hadamard test circuit is shown in Fig.~\ref{fig:Hadamard_test}. and its action is given by
\begin{equation*}
\begin{split}
\ket{0}\ket{0^a}\ket{0^s} \longrightarrow \frac{1}{2}\ket{0}\ket{0^a}\left(\ket{\bm{\varphi}}+\ket{\bm{\psi}}\right)+\frac{1}{2}\ket{0}\ket{0^a}\left(\ket{\bm{\varphi}}-\ket{\bm{\psi}}\right)+\ket{\perp''}.
\end{split}
\end{equation*}
By measuring the first two registers of the circuit shown in Fig.~\ref{fig:Hadamard_test}, we obtain the probabilities of observing the states $\ket{0}\ket{0^a}$ and $\ket{1}\ket{0^a}$, which are denoted by $p_0$ and $p_1$, respectively. We find that
\begin{equation}
p_0=\frac{\braket{\bm{\varphi}|\bm{\varphi}}+\braket{\bm{\psi}|\bm{\psi}}+2Re(\braket{\bm{\varphi}|\bm{\psi}})}{4},~~
p_1=\frac{\braket{\bm{\varphi}|\bm{\varphi}}+\braket{\bm{\psi}|\bm{\psi}}-2Re(\braket{\bm{\varphi}|\bm{\psi}})}{4}.
\end{equation}
Therefore, the real part of the inner product is
$Re(\braket{\bm{\varphi}|\bm{\psi}})=p_0-p_1$.

\begin{figure}[ht]
\[
\Qcircuit @C=1.5em @R=2em {
&\lstick{\ket{0}}       &\qw  &\gate{H} &\ctrlo{1}        &\ctrl{1}          &\gate{H}   &\meter \\
&\lstick{\ket{0^a}}    &\qw{/^a}  &\qw     &\multigate{1}{U}  &\multigate{1}{V}  &\qw        &\meter \\
&\lstick{\ket{0^s}}    &\qw{/^s}  &\qw     &\ghost{U}         &\ghost{V}         &\qw        &\qw &       \\
}
\]
\caption{ Circuit of the Hadamard test for block encodings\cite{TWY2024QCG}.}
\label{fig:Hadamard_test}
\end{figure}

\subsection{Quantum CGLS algorithm}
The classic CGLS algorithm requires repeated matrix-vector multiplications and inner product computations. These operations dominate the computational complexity, particularly for high dimensional systems. Therefore, we employ a quantum procedure to accelerate these computations.

We consider the solution of the normal equations in Eq.(\ref{eq:normal_equation}), assuming that $\|A\|\leq\alpha$. To solve this system, we employ GQSVT to obtain the solution of the rescaled system
\begin{equation}\label{eq:rescal_system}
\frac{A^{\top}A}{\alpha^2}\ket{\bm{x}}=\frac{A^{\top}}{\alpha}\frac{\ket{\bm{b}}}{\|\ket{\bm{b}}\|}
\end{equation}
in practice. Suppose the state $\frac{\ket{\bm{b}}}{\|\ket{\bm{b}}\|}$ can be efficiently prepared by an oracle $U_{b}$ and the block encoding of $A^{\top}/\alpha$ can be implemented. 

As discussed in Eq.(\ref{eq:gmp_vectors}), the vectors can be represented as
\begin{equation*}
\ket{\bm{x}_j}=\tilde{\mathcal{X}}^{\diamond}_j\left(\frac{A^{\top}}{\alpha}\right)\frac{\ket{\bm{b}}}{\|\ket{\bm{b}}\|},~
\ket{\bm{r}_j}=\tilde{\mathcal{R}}^{\diamond}_j\left(\frac{A^{\top}}{\alpha}\right)\frac{\ket{\bm{b}}}{\|\ket{\bm{b}}\|},~
\ket{\bm{p}_j}=\tilde{\mathcal{P}}^{\diamond}_j\left(\frac{A^{\top}}{\alpha}\right)\frac{\ket{\bm{b}}}{\|\ket{\bm{b}}\|}.
\end{equation*}
The absolute maximum values of $\tilde{\mathcal{X}}_j(x)$, $\tilde{\mathcal{R}}_j(x)$ and $\tilde{\mathcal{P}}_j(x)$ for $x\in[-1,1]$ are denoted as $\mathcal{X}^{\max}_{j}$, $\mathcal{R}^{\max}_{j}$ and $\mathcal{P}^{\max}_{j}$, respectively. Since the third line of CGLS algorithm requires the computation of the inner product between $\ket{\bm{p}_j}$ and $\frac{A^{\top}A}{\alpha^2}\ket{\bm{p}_j}$, we define $\mathcal{P}^{'}_{j}(t) = t^2\tilde{\mathcal{P}}_j(t)$ for notational simplicity. The maximum value of $\mathcal{P}^{'}_{j}(x)$ is denoted by $\mathcal{P}^{'\max}_{j}$. Therefore,
$$\ket{\bm{p}'_j} = \frac{A^{\top}A}{\alpha^2}\ket{\bm{p}_j} = \mathcal{P}^{'\diamond}_{j}\left(\frac{A^{\top}}{\alpha}\right)\frac{\ket{\bm{b}}}{\|\ket{\bm{b}}\|}.$$

The inner products in Algorithm \ref{algorithm.1} can be computed by employing the Hadamard test introduced in Section \ref{sect:swap-test}. As an example, we consider the computation of $\braket{\bm{p}_j|\bm{p}'_j}$. Suppose $U_{p_j}$ and $V_{p'_j}$ are GQSVT programs that encode the general matrix polynomials $\frac{1}{\mathcal{P}_j^{\max}}\tilde{\mathcal{P}}^{\diamond}_j(A^{\top}/\alpha)$ and $\frac{1}{\mathcal{P'}_j^{\max}}\mathcal{P}^{'\diamond}_j(A^{\top}/\alpha)$, respectively. That is, 
\begin{equation*}
	U_{p_j}\ket{0^a}\frac{\ket{\bm{b}}}{\|\ket{\bm{b}}\|}
	=\ket{0^a}\frac{\tilde{\mathcal{P}}^{\diamond}_j(A^{\top}/\alpha)}{\mathcal{P}_j^{\max}}\frac{\ket{\bm{b}}}{\|\ket{\bm{b}}\|}+\ket{\perp}
	=\frac{1}{\mathcal{P}_j^{\max}}\ket{0^a}\ket{\bm{p}_j}+\ket{\perp},
\end{equation*}
\begin{equation*}
	V_{p'_j}\ket{0^a}\frac{\ket{\bm{b}}}{\|\ket{\bm{b}}\|}
	=\ket{0^a}\frac{\mathcal{P}^{'\diamond}_j(A^{\top}/\alpha)}{\mathcal{P}_j^{'\max}}\frac{\ket{\bm{b}}}{\|\ket{\bm{b}}\|}+\ket{\perp}
	=\frac{1}{\mathcal{P}_j^{'\max}}\ket{0^a}\ket{\bm{p}'_j}+\ket{\perp}.
\end{equation*}
As shown in Fig.~\ref{fig:inner_product_p}, the third register of the initial state should first be acted on by $U_b$ to prepare the state $\frac{\ket{\bm{b}}}{\|\ket{\bm{b}}\|}$. By measuring the first two registers, we can estimate the value of the inner product $\braket{\bm{p}_j|\bm{p}'_j}$.

\begin{figure}[ht]
	\[
	\Qcircuit @C=1.5em @R=2em {
		&\lstick{\ket{0}}       &\qw  &\gate{H} &\ctrlo{1}        &\ctrl{1}          &\gate{H}   &\meter \\
		&\lstick{\ket{0^a}}    &\qw{/^a}  &\qw     &\multigate{1}{U_{p_j}}  &\multigate{1}{V_{p'_j}}  &\qw        &\meter \\
		&\lstick{\ket{0^s}}    &\qw{/^s}  &\gate{U_b}     &\ghost{U_{p_j}}         &\ghost{V_{p'_j}}         &\qw        &\qw &       \\
	}
	\]
	\caption{Quantum circuit for estimating the inner product $\braket{\bm{p}_j|\bm{p}'_j}$, where $U_{p_j}$ and $V_{p'_j}$ are GQSVT procedures that encode the general matrix polynomials corresponding to the states $\ket{\bm{p}_j}$ and $\ket{\bm{p}'_j}$, respectively.}
	\label{fig:inner_product_p}
\end{figure}

We propose a hybrid classical-quantum algorithm, called quantum CGLS, as shown in Algorithm \ref{algorithm.2}. Step 1 employs the Hadamard test to estimate the inner products, which serves as the starting point of the iteration. Steps 2-19 form the main iterative framework. Specifically, steps 10 and 18 invoke the Hadamard test, while the remaining steps are classical updates of the polynomial coefficients. The states involved in the inner products are prepared by GQSVT, while the corresponding phase factors are computed classically. The stopping criterion in step 11 is proven in subsection \ref{sect:circuit_depth}. After the iteration is completed, the approximated solution of the system in Eq.~(\ref{eq:rescal_system}) is prepared by GQSVT in steps 20 and 21.

\begin{algorithm}[ht]
\caption{Quantum CGLS algorithm }
\label{algorithm.2}
\textbf{Input:} $\ket{\bm{x}_0}=0$, 
an error tolerance $\varepsilon$.\\
\textbf{Output:} a state $\ket{\bm{x}}$ approximating the solution of $A^{\top}A\ket{\bm{x}}=A^{\top}\ket{\bm{b}}$.\\
\textbf{Procedure:}
\begin{algorithmic}[1]
\STATE Estimate  $\braket{\bm{r}_0|\bm{r_0}}$ and $\braket{\bm{p}_0|\bm{p}'_0}$ using  the Hadamard test, 
where $\ket{\bm{r}_0}=\ket{\bm{p}_0} = \frac{A^{\top}}{\alpha}\frac{\ket{\bm{b}}}{\|\ket{\bm{b}}\|} $,
$\mathcal{P'}_0(t) = t^3$ and $\ket{\bm{p_0}'} = \mathcal{P}^{'\diamond}_0\left(\frac{A^{\top}}{\alpha}\right)\frac{\ket{\bm{b}}}{\|\ket{\bm{b}}\|} $.

\FOR {$j=0,1,...,m$ }
   \STATE Perform the following steps $4-9$ on a classical computer.
   \STATE  $\alpha_j=\frac{\braket{\bm{r}_j|\bm{r}_j}}{\braket{\bm{p}_j|\bm{p}'_j}}$
   \FOR {$l=0,1,...,j+1$} 
      \STATE $\chi_{l}^{(j+1)}=\chi_l^{(j)}+\alpha_j\rho_l^{(j)}$
      \STATE $\gamma_{l}^{(j+1)}=\gamma_l^{(j)}-\alpha_j\rho_{l-1}^{(j)}$
   \ENDFOR
   \STATE  Compute the maximum value $\mathcal{R}^{\max}_{j+1}$ of $\tilde{\mathcal{R}}_{j+1}(t)$ and the phase factors for the normalized polynomial.
   \STATE  Estimate  $\frac{\braket{\bm{r}_{j+1}|\bm{r}_{j+1}}}{(\mathcal{R}^{\max}_{j+1})^2}$ using the Hadamard test, where $\ket{\bm{r}_{j+1}}$ is obtained by GQSVT according to the  polynomial $\tilde{\mathcal{R}}_{j+1}(t)/\mathcal{R}^{\max}_{j+1}$ and its corresponding phase factors.
   \STATE \textbf{if} {$\|\bm{r}_{j+1}\|\leq \frac{\varepsilon\|A\|^2}{\alpha\|\ket{\bm{b}}\|\kappa(A)^2}$} \textbf{then} break
   \STATE Perform the following steps $13-17$ on a classical computer.
   \STATE $\beta_j=\frac{\braket{\bm{r}_{j+1}|\bm{r}_{j+1}}}{\braket{\bm{r}_j|\bm{r}}_j}$
   \FOR{$l=0,1,...,j+1$} 
      \STATE $\rho_{l}^{(j+1)}=\gamma_l^{(j)}+\beta_j\rho_l^{(j)}$
   \ENDFOR
   \STATE  Compute the maximum value $\mathcal{P}^{\max}_{j+1}$ and $\mathcal{P}^{'\max}_{j+1}$ of $\tilde{\mathcal{P}}_{j+1}(t)$ and $\mathcal{P}'_{j+1}(t)$. Then determine the phase factors of these normalized polynomials.
   \STATE Estimate $\frac{\braket{\bm{p}_{j+1}|\bm{p}'_{j+1}}}{\mathcal{P}^{\max}_{j+1}\mathcal{P}^{'\max}_{j+1}}$ using the Hadamard test, where $\ket{\bm{p}_{j+1}}$ and $\ket{\bm{p}'_{j+1}}$ are implemented via GQSVT according to the polynomials $\tilde{\mathcal{P}}_{j+1}(t)/\mathcal{P}^{\max}_{j+1}$ and $\mathcal{P}'_{j+1}(t)/\mathcal{P}^{'\max}_{j+1}$ and their phase factors.
\ENDFOR
\STATE  Compute the maximum value $\mathcal{X}^{\max}_{m+1}$ of  $\tilde{\mathcal{X}}_{m+1}$ and the phase factors for the normalized polynomial.
\STATE Prepare the state proportional to $\frac{1}{\mathcal{X}^{\max}_{m+1}}\ket{\bm{x}_{m+1}}$ by GQSVT according to the polynomial $\tilde{\mathcal{X}}_{m+1}/\mathcal{X}^{\max}_{m+1}$ and its corresponding phase factors.
\end{algorithmic}	
\end{algorithm}

\subsection{The maximum circuit depth}\label{sect:circuit_depth}
It is well known that the convergence of the CG algorithm for the SPD system $\tilde{A}\tilde{\bm{x}} = \tilde{\bm{b}}$ satisfies
\begin{equation}\label{eq:CGconvergence}
	\left\|\tilde{\bm{x}}-\tilde{\bm{x}}_m\right\|_{\tilde{A}}\leq 2\left(\frac{\sqrt{\kappa}-1}{\sqrt{\kappa}+1}\right)^m
	\left\|\tilde{\bm{x}}-\tilde{\bm{x}}_0\right\|_{\tilde{A}},
\end{equation}
where $\|\tilde{\bm{x}}\|_{\tilde{A}} = \sqrt{\tilde{\bm{x}}^{\top}\tilde{A}\tilde{\bm{x}}}$, $\tilde{A}$ denotes the SPD coefficient matrix and $\kappa$ is the condition number of $\tilde{A}$\cite{Saad2003}. We now employ it to calculate the maximum circuit depth of our algorithm.

Solving the rescaled system in Eq.~(\ref{eq:rescal_system}) yields an approximate solution vector $\ket{\bm{x}_m}$, which approximates $\frac{\alpha}{\|\ket{\bm{b}}\|}\bm{x}$, where $\bm{x}$ (without Dirac notation) is the exact solution of the normal equations $A^{\top}A\bm{x} = \bm{b}$. 
Applying the CG convergence formula in Eq.~(\ref{eq:CGconvergence}) and let $\tilde{A} = \frac{A^\top A}{\alpha^2}$, we have
\begin{equation}\label{eq:CGLSconvergence}
\left\|\frac{\alpha}{\|\ket{\bm{b}}\|}\bm{x}-\ket{\bm{x}_m}\right\|_{\frac{A^\top A}{\alpha^2}}\leq 2\left(\frac{\sqrt{\kappa(A^\top A)}-1}{\sqrt{\kappa(A^\top A)}+1}\right)^m
\left\|\frac{\alpha}{\|\ket{\bm{b}}\|}\bm{x}-\ket{\bm{x}_0}\right\|_{\frac{A^\top A}{\alpha^2}}.
\end{equation}
Note that $\ket{\bm{x}_0}=\bm{0}$, $\kappa(A^\top A)=\kappa(A)^2$, and
\begin{equation}\label{eq:norm_scaling}
\left\|\bm{x}\right\|_{\frac{A^\top A}{\alpha^2}}
= \left\|\frac{A}{\alpha}\bm{x}\right\|
\leq\left\|\bm{x}\right\|
\leq \frac{\kappa(A^\top A)\|A^{\top}\ket{\bm{b}}\|}{\|A^\top A\|} = \frac{\kappa(A)^2\|A^{\top}\ket{\bm{b}}\|}{\|A\|^2}.
\end{equation}
Then, the Eq.(\ref{eq:CGLSconvergence}) becomes 
\begin{equation*}
\left\|\bm{x}-\frac{\|\ket{\bm{b}}\|}{\alpha}\ket{\bm{x}_m}\right\|_{\frac{A^\top A}{\alpha^2}}\leq 2\left(\frac{\kappa(A)-1}{\kappa(A)+1}\right)^m \frac{\kappa(A)^2\|A^{\top}\ket{\bm{b}}\|}{\|A\|^2}.
\end{equation*}
Thus, the number of iterations can be calculated as
\begin{equation*}
	m=O\left(\kappa(A)\log\left(\frac{\kappa(A)^2\left\|A^{\top}\ket{\bm{b}}\right\|}{\varepsilon\|A\|^2}\right)\right)
\end{equation*}
for the pre-specified error $\varepsilon$. Moreover, when performing the $m$-th iteration, the GQSVT with a polynomial of degree at most $m+1$ is implemented, resulting in a maximum circuit depth of
\begin{equation}
 2m+3=O\left(\kappa(A)\log\left(\frac{\kappa(A)^2\left\|A^{\top}\ket{\bm{b}}\right\|}{\varepsilon\|A\|^2}\right)\right).
\end{equation}

In addition, the stopping criterion can be formulated in terms of the residual norm. For the purpose of 
\begin{equation}
\begin{split}
\left\|\bm{x}-\frac{\|\ket{\bm{b}}\|}{\alpha}\ket{\bm{x}_m}\right\|
& \leq \left\|(A^\top A)^{-1}\right\|\cdot\left\|A^{\top}\ket{\bm{b}}-\frac{\|\ket{\bm{b}}\|}{\alpha}(A^\top A)\ket{\bm{x}_m}\right\|\\
& = \alpha\left\|\ket{\bm{b}}\right\|\left\|(A^\top A)^{-1}\right\|\cdot\left\|\frac{A^{\top}}{\alpha}\frac{\ket{\bm{b}}}{\|\ket{\bm{b}}\|}-\frac{A^\top A}{\alpha^2}\ket{\bm{x}_m}\right\|\\
& = \alpha\left\|\ket{\bm{b}}\right\|\left\|(A^\top A)^{-1}\right\|\cdot\left\|\bm{r}_m\right\|\leq\varepsilon,
\end{split}
\end{equation}
the stop condition is set to $\left\|\bm{r}_m\right\|
\leq
\frac{\varepsilon}{\alpha\|\ket{\bm{b}}\|\|(A^\top A)^{-1}\|} 
= \frac{\varepsilon\|A\|^2}{\alpha\|\ket{\bm{b}}\|\kappa(A)^2}$. Here,  $\|\cdot\|_{\frac{A^\top A}{\alpha^2}}$ is replaced by the $\ell^2$-norm $\|\cdot\|$ because of Eq.(\ref{eq:norm_scaling}).

{\bfseries Remark 2.} Solving the linear system (\ref{eq:rescal_system}) directly using the QSVT method yields a maximum circuit depth of $O(\kappa(A)^2\log(\kappa(A)^2/\varepsilon))$. Consequently, our approach achieves a quadratic speedup in the maximum circuit depth with respect to the condition number.

\section{Conclusion}
In this work, we present the GQSVT method and apply it to the classical CGLS algorithm, thereby proposing a quantum CGLS algorithm. As an extension of QSP, the GQSP removes additional constraints on implementable polynomial families, retaining only the necessary constraint of $|p|\leq 1$, but GQSP only works for unitary matrices. We further generalize GQSP by employing four types of controlled unitary as signal operators, combining it with qubitization and singular value decomposition, thereby realizing the block encoding of indefinite parity polynomials for general operations. Additionally, since CGLS is a Krylov subspace method, when solving the normal equations $A^{\top}A\ket{\bm{x}}=A^{\top}\ket{\bm{b}}$ using the conjugate gradient method, the iterative vectors are expressed as matrix functions of $A^{\top}A$ acting on the right-hand side $A^{\top}\ket{\bm{b}}$. They can then be reformulated as generalized matrix functions of $A$ acting on $\ket{\bm{b}}$, enabling us to utilize GQSVT. Finally, we analyze the maximum circuit depth and show that our method achieves a quadratic speedup with respect to the matrix condition number.

\bibliographystyle{elsarticle-num}
\bibliography{reference_GQSVT}

\end{document}